\documentclass[compress,final,3p,times,11pt]{elsarticle}
\usepackage{amsfonts}
\usepackage{dsfont}
\usepackage{amssymb}
\usepackage{amsmath}
\usepackage{color,epstopdf}
\usepackage{mathrsfs}
\usepackage{graphicx,color,epstopdf}
\usepackage{float}
\usepackage{caption}
\usepackage{subcaption}
\usepackage{float}
\usepackage{subcaption}

\newtheorem{theorem}{Theorem}
\newtheorem{remark}{Remark}

\newtheorem{definition}{Definition}

\journal{Nonlinear Dynamics}

\begin{document}
\begin{frontmatter}

\title{Most Probable Dynamics of the Single-Species with Allee Effect under Jump-diffusion Noise}

\author[addr1,addr2]{Almaz Tesfay}\ead{almaz.tesfay@mu.edu.et}
\author[addr1]{Shenglan Yuan\corref{cor1}}\ead{shenglanyuan@hust.edu.cn}\cortext[cor1]{Corresponding author}
\author[addr1,addr2,addr4]{Daniel Tesfay}\ead{daniel.tesfay@mu.edu.et}
\author[addr3]{James Brannan}\ead{jrbrn@clemson.edu}

\address[addr1]{Center for Mathematical Sciences, Huazhong
University of Science and Technology, Wuhan 430074, China}
\address[addr2]{Department of Mathematics, Mekelle University, P.O.Box 231, Mekelle, Ethiopia}
\address[addr3]{Department of Mathematical Sciences, Clemson University, Clemson, South Carolina 29634, USA}
\address[addr4]{ Department of Rural Development and Agricultural Economics, University of Rwanda,  P.O. Box 210 , Musanze, Rwanda }



\begin{abstract}
We investigate the most probable phase portrait (MPPP) of a stochastic single-species model with the Allee effect using the non-local Fokker-Planck equation. This stochastic model is driven by non-Gaussian as well as Gaussian noise, and it has three fixed points. One of them is the unstable state which lies between the two stable equilibria. We focus on the transition pathways from the extinction state to the upper fixed stable state for the transcription factor activator in a single-species model. This helps us to study the biological behavior of species.
 The most probable path is obtained from the solution of the non-local Fokker-Planck equation corresponding to the population system of the single-species model, and the corresponding maximum possible stable equilibrium state is determined. We also obtain the Onsager-Machlup (OM) function for the stochastic model and solve the corresponding most probable paths. The numerical simulation shows that: (i) When non-Gaussian noise is presented in the system, the maximum of the stationary density function is located at the most probable stable equilibrium state; (ii) If the initial value increases from extinction state to the upper stable state, the most probable trajectory goes to the maximal likely equilibrium state, in our case it lies between 9 and 10; (iii) The most probable paths increase to stable state quickly, then maintain a nearly constant level, and approach to the upper stable equilibrium state as time goes on. These numerical experiment findings accelerate growth for further experimental study, in order to achieve good knowledge about dynamical systems in biology.

\end{abstract}




\end{frontmatter}
\section{Introduction}
Single-species dynamics is one of the core research areas in theoretical ecology. Research about single-species dynamics enables the researcher to find out the conditions of extinction and persistence of the species.  The strong motivation for the researchers to develop mathematical models is to understand the cause of cycles, like populations \cite{Yz1}.

Population modeling is very important for species management, for example, in developing recovery plans for species threatened by extinction, managing fisheries for the highest possible sustainable yield, and trying to contain or prevent the spread of invasive species \cite{r0,b1,b0}.

In literature, one can find several models of the dynamical single-species growth system. Gompertz growth model \cite{gq}, Verhulst growth model with or without Allee effect \cite{lg}, power law growth model \cite{pl}, {the interconnections between deterministic and stochastic system \cite{new}}, Gilpin--Ayala  model \cite{ga} are only a few to mention.  In this study, we compute a single-species model focusing on the Verhulst growth model with the Allee effect developed by Y. Jin \cite{7}.

The dynamics of biological phenomena, particularly that of populations of living beings, besides some clear trends, are frequently influenced by unpredictable components due to the complexity and variability of environmental conditions \cite{r3}. Extensive researchers in modeling and analysis of random fluctuations \cite{r41,r4,TSB} in biological dynamical systems have been ongoing for a long time now. The studies of events in the population such as persistence stationary distribution, and extinction in stochastic single-species models become an interesting and important research field. One of the hot issues in population dynamics is developing sufficient conditions for the persistence of biological species as mentioned in \cite{r5,r7,r6} and the references therein.

The population may be affected by sudden environmental noises \cite{rr0,YLZ}. For example, severe acute respiratory syndrome (SARS),  human immunodeficiency virus (HIV), the smoking habit \cite{smoke}, and the recent COVID-19 \cite{ab1}, earthquakes \cite{YBD}, temperature \cite{YZD}, and hurricanes \cite{YW23}. These sudden environmental perturbations may bring substantial social and economic losses. Stochastic single-species model perturbed by Brownian motion has been researched extensively by many authors \cite{rr1,rr3,rr2,amu1,ab}. However, stochastic extension of population process driven by Gaussian noise cannot explain the aforementioned random and intermittent environmental perturbations. Introducing a L\'{e}vy process into the underlying population dynamics would explain the impact of these random jumps  \cite{l1}. There have been a few studies that investigated dynamical systems where the noise source is a L\'{e}vy process \cite{l2,l3,YD19,ga}. Implying the L\'evy noise into the biological system to simulate the effect caused by the external environment is more effective and nearer to reality than the Gaussian noise. The investigation of the single-species model is still in its infancy even though noisy fluctuations naturally portray random intermittent jumps.  L\'{e}vy noise is widely applied in studying natural and man-made phenomena in science, among which we mention biology \cite{c1},  physics \cite{s1,YB22},  and economics \cite{s2,s3}.

 Under this research, we consider the population dynamics of a single-species growth model with Allee effect perturbed by stable L\'evy fluctuations.  We also analyze the influence of L\'evy  noise fluctuation on the system (\ref{401}).  Investigating the impact of noisy fluctuations acts a pivotal part in demonstrating the intricate interactions between the single-species models and their complex surroundings. We study how Allee effects and stochasticity combine to affect population persistence in here. To find the numerical solutions for the Fokker-Planck equation  determined by non-local differential equation with symmetric $\alpha$-stable L\'evy motion, we  apply a finite difference method probed  by Gao et al. \cite{Gao}.

The most probable phase portrait was first proposed by Duan \cite[Section 5.3.3]{1}. Cheng et al. \cite{2} obtained the analytical results of the MPPP and showed that the MPPP can give  useful information about the propagation of stochastic dynamics in the one-dimensional model. Wang et al. \cite{3} studied the stochastic bifurcation by using the qualitative changes of the MPPP to a stochastic system driven by multiplicative stable L\'evy noise. In Ref. \cite{4}, the  scholars  investigated  the most probable trajectories of the tumor growth system with immune surveillance under correlated Gaussian noises, and  derived analytical solution of the most probable steady state by utilizing the extremum theory with the local Fokker-Planck equation (FPE) in the system.  A function which summarizes about the behavior of the  dynamics of a continuous stochastic process was defined as the Onsager-Machlup function \cite{15}.
The Onsager-Machlup function for stochastic models driven by both non-Gaussian and Gaussian noises was established in \cite{5}. The authors also examined  the corresponding MPPP of the stochastic dynamical  systems. Cheng et al. \cite{6} focused on the impact of Gaussian noise and jump stable
L\'evy noise in a genetic regulatory system, and they minimized the OM action functional for the stochastic dynamics driven by Gaussian and obtained the most probable transition pathway.  This inspired us to study the MPPP of the single-species model.

Therefore, our goal is about to investigate how the most probable trajectories escape from the single-species state to the extinction state more quickly.

Consider the following stochastic single-species growth model with Allee effect:
{\small\begin{equation}\label{401}
dX_t=X_{t-}\left[\left(s-\gamma_2\,X_{t-}-\frac{\gamma_3}{\gamma_3\,\gamma_4\,X_{t-}+1}\right)dt + \lambda dB_t+ \int_\mathbb{Y}\epsilon(y)\tilde{N}(dt,dy)\right],
\end{equation}}
 for $t\geq 0$ and $X_0=x_0$, where $X_{t-}$ is the left limit of the population size $X_t$.
\begin{table}[h!]
\centering
\parbox{0.8\textwidth}{
\caption{Biological meaning of the parameters and variables in the single-species model.}
\begin{tabular}[t]{|l |l|l| p{1.1cm}}

\hline
     Parameter &Definition \\ \hline \hline
     $s$&  The growth rate\\ \hline
     $\gamma_2$& Intraspecific competition rate\\ \hline
     $\gamma_3$& The attack rate\\ \hline
     $\gamma_4$& Represents the handling time of predator\\ \hline
     $M=s/\gamma_2$& The carrying capacity\\ \hline
     $t$ & Time\\ \hline
\end{tabular}
}
\end{table}\label{P}

Stochastic force $\tilde{N}(dt,dy)=N(dt,dy)-\nu_{\alpha}(dy)dt$ is a compensated Poisson random measure with associated Poisson random measure $N(dt,dy)$ and intensity measure $\nu_{\alpha}(dy)dt$, in which $\nu_{\alpha}(dy)$ is L\'evy measure on a measurable subset $\mathbb{Y}$ of $(0,\infty)$ with $\nu_{\alpha}(\mathbb{Y})<\infty$; see \cite{YSD}.

The following restriction on system (\ref{401}) is natural for biological meaning:
\begin{equation*}
1+\epsilon(y)>0,\quad y\in\mathbb{Y}.
\end{equation*}
When $\epsilon(y)>0$, the perturbation stands for the increasing species, e.g. planting, while $\epsilon(y)<0$ represents that the species is decreasing, e.g. harvesting and epidemics.

The main aim of this study is to investigate stochastic dynamics of single-species biological populations in random environments. We model the evolution of these populations with
first order ordinary autonomous differential equations bringing in the coefficients and inputs which are stochastic processes. The two stochastic processes germane to this study are Brownian motion and L\'evy process. Brownian motion describes random fluctuations that are continuous in time; see Subsection 2.1. L\'evy process, of which Brownian motion is a special case, is used to model random fluctuations which may have discontinuities or jumps; see Subsection 2.2.

Here, we develop the stochastic single-species model with the Allee effect influenced by Gaussian and non-Gaussian noises. Firstly, we review the deterministic model, calculate its equilibrium solutions and describe the behavior of the fixed points. Secondly, we get the highest possible paths, and the corresponding maximum possible stable states attracting the nearby maximum possible paths of the stochastic system (\ref{401}). We do this by finding the stationary density function which is the solution of the non-local FPE. To solve the non-local partial differential equation, we use the finite difference method proposed in \cite{Gao}. This method helps us explore some dynamical behaviors of the single-species system under the impact of non-Gaussian L\'evy noise.

This study is organized as follows. In the second section, we recall the definitions of the one-dimensional Brownian motion $B_t$ and  symmetric $\alpha$-stable L\'evy motion $L_{t}^{\alpha}$. In the third section, we discuss the formulation and  analysis of the deterministic model (\ref{42}) of the single-species with Allee effect . In the fourth section, we explain the analysis of the stochastic single-species model (\ref{401}) with Allee effect. We also review the definition of the Onsager-Machlup function and most probable phase portraits in the subsections \ref{OM} and \ref{Mppp}, respectively. The numerical results and the biological implication of our experimental findings are presented in the fifth section. We conclude our research by  giving a brief summary in the last section.
\section{Preliminaries}\label{Sec2}
Under this section, we define the one-dimensional Brownian motion starting at time $t=0$ as a process $B_t$ and $\alpha$-stable L\'evy motion $L_t^{\alpha}$, which constitute a class of stochastic processes that have independent and stationary increments as defined below.  Throughout this study, we denote $\mathbb{R}^+=(0,\infty)$, \, $\mathbb{R} = (-\infty, \infty)$, and $X_t \in \mathbb{R}^+$, for $t\geq 0$.
\subsection{\textbf{{Brownian motion}}}\label{Sec22}
Brownian motion $B_t$ (also called Wiener process) is a one-dimensional stochastic process defined on complete probability space $(\Omega, \mathcal{F},\mathcal{F}_t,\mathbb{P})$, which has independent and stationary increments \cite{mao,11}. Brownian motion $B_t$ satisfies the following conditions: \\(i) $B_t$ has continuous paths, and its paths are nowhere differentiable almost surely;\\(ii) $B_t$ has stationary increments, i.e.,  $B_t-B_s$ is normally distributed with mean 0 and variance $t-s$ for any $0\leq s\leq t$;\\(iii) The process starts at the origin, i.e., $B_0=0$ almost surely;\\
(iv) $B_t$ has independent increments, i.e., $B_t-B_s$ is independent of the past for $s<t$.

\subsection{\textbf{The $\alpha$-Stable L\'evy motion}}\label{Sec23}
L\'evy motions $L_t$ are a class of non-Gaussian stochastic processes. A  L\'evy motion $L_t$ having values in $\mathbb{R}$ is determined by a drift coefficient $\hat{b}\in \mathbb{R}$, $\hat{Q}\geq0$ and a Borel measure $\nu$ defined on ${\mathbb{R}}\setminus{\{0\}}$. The triplet $(\hat{b},\hat{Q},\nu)$ is the so-called generating triplet of L\'evy motion $L_t$.  A L\'evy motion can be written as linear combination of time $t$, a Brownian motion and a pure jumping
process \cite{12,11}, i.e., $L_t$ can be expressed as
\begin{equation*}\label{Ito-Dec}
 L_t=\hat{b}t+ B_{\hat{Q}}(t)+\int_{|y|<1}y\tilde{N}(t,dy)+\int_{|y|\geq 1}yN(t,dy),
\end{equation*}
where $N(t,dy)$ is the independent  Poisson random measure on $\mathbb{R}^+\times{{\mathbb{R}}\setminus{\{0\}}}$, $\tilde{N}(t,dy)=N(t,dy)-\nu(dy)dt$ is the compensated Poisson random measure, $\nu(S)=\mathbb{E}(N(1,S))$ is the jump measure, and $B_{\hat{Q}}(t)$ is the independent Brownian motion.

The L\'evy-Khinchin formula for L\'evy motion has a specific form of its characteristic function.
\begin{equation*}
\mathbb{E}[e^{(i\,\xi{L_t})}] = e^{t\phi(\xi)},\quad 0\leq t<\infty,
\end{equation*}
where
$$\phi(\xi)=i\,\xi\hat{b}-\frac{\hat{Q}}{2}\,\xi^2 + \int_{{\mathbb{R}}\setminus{\{0\}}}(e^{i\,\xi z}-1-i\,\xi z\mathds{1}{_{|z|<1}})\nu(dz), \quad \xi \in \mathbb{R}.$$

A stable distribution $S_{\alpha}(\theta,\beta,\gamma)$ is the distribution for a stable random variable, where the stability index $\alpha \in (0,2)$, the skewness $\beta \in (0,\infty)$, the shift $\gamma \in (-\infty,\infty)$, and scale index $\theta \geq0$. A $\alpha$-stable L\'evy motion $L_{t}^{\alpha}$ \cite{9,8,10} is a non-Gaussian stochastic process satisfying\\
(i) the random variables $L_{t_{i+1}}^{\alpha}-L_{t_{i}}^{\alpha}$ are independent for $0\leq t_{1}<t_{2}<\cdot\cdot\cdot<t_{i-1}<t_{i}<t_{i+1}<\infty,$ and for each $i=1,2,\cdot\cdot\cdot$;\\
(ii) $L_{t}^{\alpha}$ has stochastically continuous sample paths, i.e., for $0 \leq s\leq t$ and $\delta>0$, the probability $\mathbb{P}(|L_{t}^{\alpha}-L_{s}^{\alpha}|>\delta)$ approaches to zero as $t\rightarrow s$;\\
(iii) $L_{0}^{\alpha}=0$, almost surely;\\
(iv) $L_{t}^{\alpha}-L_{s}^{\alpha}$ and $L_{t-s}^{\alpha}$ have the same distribution $S_{\alpha}((t-s)^{1/{\alpha}},0,0)$.

In the case of a one-dimensional isotropic $\alpha $-stable L\'evy motion, the L\'evy triplet has the drift factor $\hat{b}=0$ and the diffusion coefficient $\hat{Q}=0$.
In this study, we focus on jump process with a specific size in generating triplet $(0,0,\nu_{\alpha})$ for the random distribution $S_{\alpha}$ which can be defined by $\Delta L_{t}^{\alpha}=L_{t}^{\alpha}-L_{t-}^{\alpha}<\infty, \,\,t \geq 0,$ where $L_{t-}^{\alpha}$ is the left limit of the $\alpha$-stable L\'evy motion in $\mathbb{R}$ at any time $t$.  Here
$\nu_{\alpha}(dz)=c({\alpha})\frac{1}{|z|^{1+\alpha}}dz$ is L\'evy measure
with $ c_{\alpha} = \alpha\frac{\Gamma(\frac{1+\alpha}{2})}{{2^{1-\alpha}\pi^{\frac{1}{2}}}\Gamma{(1-\frac{\alpha}{2})}}$, and $\Gamma$ is the Gamma function.
\begin{remark}
A special case of $\alpha$-stable L\'evy motion is Brownian motion when $\alpha$ = 2. Poisson process, $\alpha$-stable process, compound Poisson process, etc. are also examples of L\'evy processes \cite{11}.
\end{remark}
\section{Dynamical analysis of the deterministic model}\label{model}
The deterministic form of the nonlinear model (\ref{401}) without noise is given as
\begin{equation}\label{42}
\frac{dX_t}{dt}=X_{t}\left(s-\gamma_2\,X_{t}-\frac{\gamma_3}{\gamma_3\,\gamma_4\,X_{t}+1}\right)=:F(X_t), \quad t\geq 0,\quad X_0=x_0.
\end{equation}
This system can be written as $\frac{dX}{dt}=-\frac{dU(X)}{dX}$, where $U(X)$ is the potential function given by
\begin{align*}
    U(X):=-\frac{s\,X^2}{2}+\frac{\gamma_2\,X^3}{3}+\frac{\gamma_3}{(\gamma_3\,\gamma_4)^2}[\gamma_3\,\gamma_4\,X+1-\ln(\gamma_3\,\gamma_4\,X+1)].
\end{align*}
The single-species model (\ref{42}) with Allee effect  has equilibrium points $X_1=0$ and
\begin{align*}
  X_{2,3}&=\frac{(s\gamma_3\gamma_4-\gamma_2)\pm \sqrt{(s\gamma_3\gamma_4-\gamma_2)^{2}-4\gamma_2\gamma_3\gamma_4(\gamma_3-s)}}{2\gamma_2\gamma_3\gamma_4}\nonumber\\
  &=  \frac{(s\gamma_3\gamma_4-\gamma_2)\pm (s\gamma_3\gamma_4-\gamma_2)\,\sqrt{1-\beta}}{2\gamma_2\gamma_3\gamma_4}\nonumber\\
 &=\frac{(s\gamma_3\gamma_4-\gamma_2)\,\Big(1\pm \,\sqrt{1-\beta}\,\,\Big)}{2\gamma_2\gamma_3\gamma_4},
\end{align*}
where  $\beta=\frac{4\gamma_2\gamma_3\gamma_4}{(s\gamma_3\gamma_4-\gamma_2)^{2}}\,(\gamma_3-s).$ If $\beta <1$, then the equilibrium states of system (\ref{42}) are
\begin{align*}
  &X_1=0,\quad \text{an extinction equilibrium;}  \nonumber\\
  &X_2=\frac{(s\gamma_3\gamma_4-\gamma_2) \,\Big(1-\sqrt{1-\beta}\Big)}{2\gamma_2\gamma_3\gamma_4},\quad \text{a lower unstable equilibrium;}  \nonumber\\
  &X_3=\frac{(s\gamma_3\gamma_4-\gamma_2) \,\Big(1+\sqrt{1-\beta}\Big)}{2\gamma_2\gamma_3\gamma_4}, \quad \text{an upper stable equilibrium}.
\end{align*}
If $\beta=1$, the single-species deterministic model (\ref{42}) has only two equilibria states:
\begin{align*}
 \text{stable state}~~ X_1=0,\quad    \text{and} \quad  \text{unstable state}~~ X_4=\frac{s\gamma_3\gamma_4-\gamma_2}{2\gamma_2\gamma_3\gamma_4}.
\end{align*}

\begin{figure}[H]
	\begin{minipage}{0.48\linewidth}
		\leftline{(a)}
		\centerline{\includegraphics[height = 6cm, width = 8.6cm]{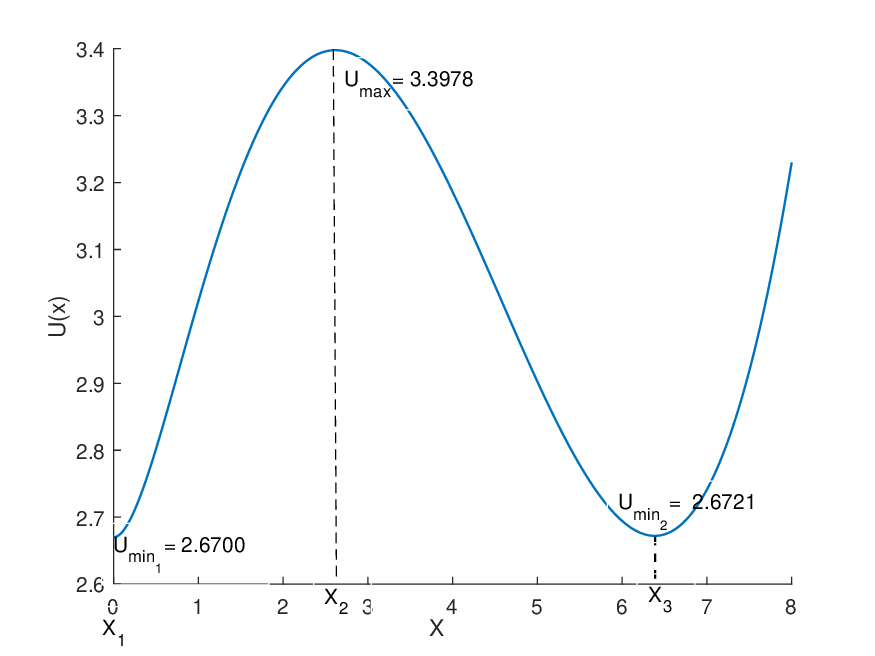}}
	\end{minipage}
	\hfill
	\begin{minipage}{0.48\linewidth}
		\leftline{(b)}
		\centerline{\includegraphics[height = 6cm, width = 8.6cm]{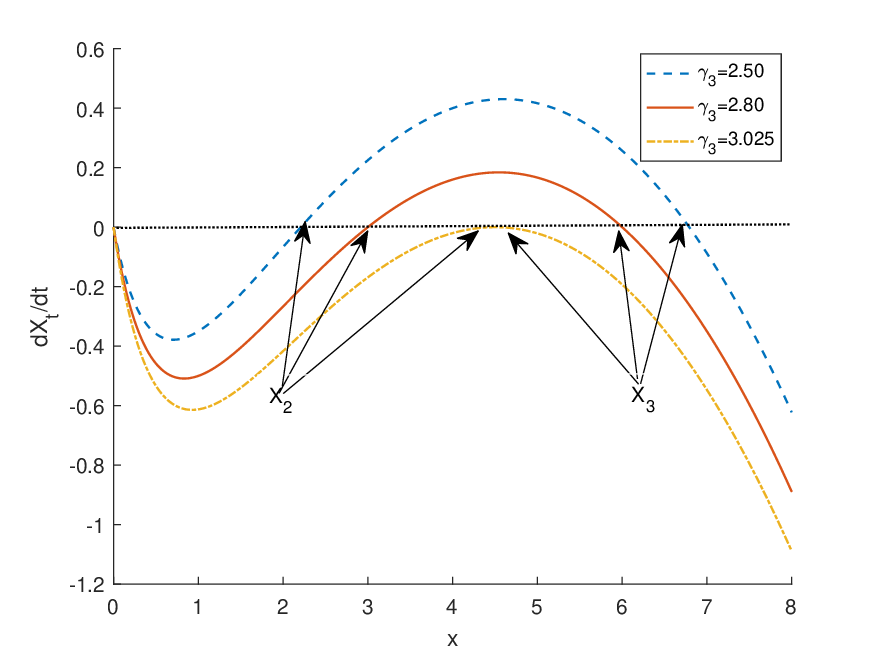}}
	\end{minipage}
\caption{(a) Numerical simulation of the bistable potential function $U(X)$ of the nonlinear model (\ref{42}). Dashed black lines: local unstable and stable equilibria at $X_2=2.6159$ and $X_3=6.3841$, respectively. (b) The phaselines of the single-species model (\ref{42}). Parameters $s=1,\,\gamma_2=0.1,\,\gamma_3=2.67,\,\gamma_4=\frac{1}{\gamma_3}$ in the graph of $\frac{dX}{dt}$.}\label{de-pot-bif}
\end{figure}

The derivative of $F(X)$ is
\begin{equation*}
    s-2\gamma_2X-\frac{\gamma_3}{(\gamma_3\,\gamma_4\,X+1)^2}.
\end{equation*}
For simplicity and  convenience of discussion, we choose the parameters $\gamma_3\,\gamma_4=1$, $s=1$, $0< \gamma_2< 1$, $0< \gamma_3< \frac{(1+\gamma_2)^2}{4\,\gamma_2}$, therefore $\beta=\frac{4\gamma_2(\gamma_3-1)}{(1-\gamma_2)^{2}},$\,\, and $X_4=\frac{1-\gamma_2}{2\,\gamma_2}$.
 For $\beta <1$, the extinction state $X_1=0$, and the equilibrium solution $X_3$ are stable, but $X_2$ is unstable. Fig. \ref{de-pot-bif}(b) shows that when the value of attack rate $\gamma_3$ increases, the unstable state $X_2$ and stable state $X_3$ get more close to each other, then become one solution and finally disappear indicating the occurrence  of the saddle-node bifurcation.

The critical value of attack rate $\gamma_c=2.67$ of the deterministic single-species system (\ref{42}) with Allee effect  is obtained by solving the equation $U(X_1)=U(X_3)$. This value is an indication to the transition phenomena between {the unstable and stable state for deterministic single-species growth model. The steady state (extinction state) $X_1$ is stable if $\gamma_3>\gamma_c$, and the steady state $X_3$ exhibits the stability property for  $\gamma_3 <\gamma_c$.}

\section{Dynamical analysis of the stochastic system}\label{stoch}
In this section, we discuss the behavior of the solution of the stochastic system (\ref{401}). Firstly, we  recall the definition of the Onsager-Machlup function for the stochastic differential equation driven by jump noise. This helps to measure OM induced by the jump process. Secondly, we  examine the corresponding most probable paths.  {Finally, we present the numerical experiment findings using finite difference method \cite{Gao}}. Hence, the numerical solution of the stochastic model  provides useful information for understanding the dynamical behavior of the system (\ref{401}).

\subsection{Onsager-Machlup  functional}\label{OM}
{The Onsager-Machlup functional defines a probability density for a stochastic process in which the probability density is estimated implicitly. It can be used for purposes of reweighting and sampling trajectories, as well as determining the most probable trajectory based on variational arguments. The most probable transition pathway can be obtained by minimizing the Onsager-Machlup function. The whole procedure enables us to detect the dynamics of the most probable path \cite{HCYD}.}

\begin{figure}[H]
	\begin{minipage}{0.48\linewidth}
		\leftline{(a)}
		\centerline{\includegraphics[height = 6cm, width = 8.6cm]{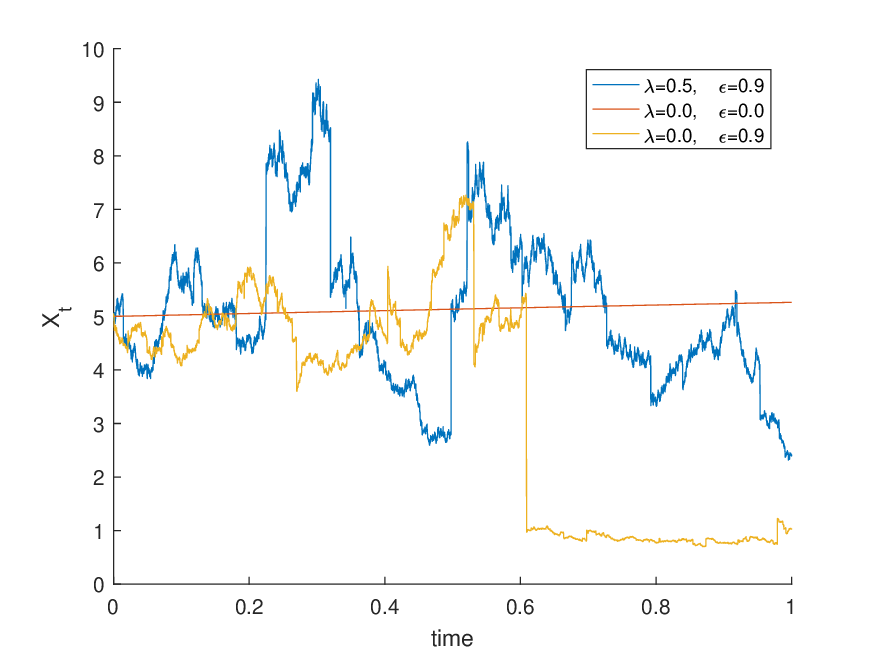}}
	\end{minipage}
	\hfill
	\begin{minipage}{0.48\linewidth}
		\leftline{(b)}
		\centerline{\includegraphics[height = 6cm, width = 8.6cm]{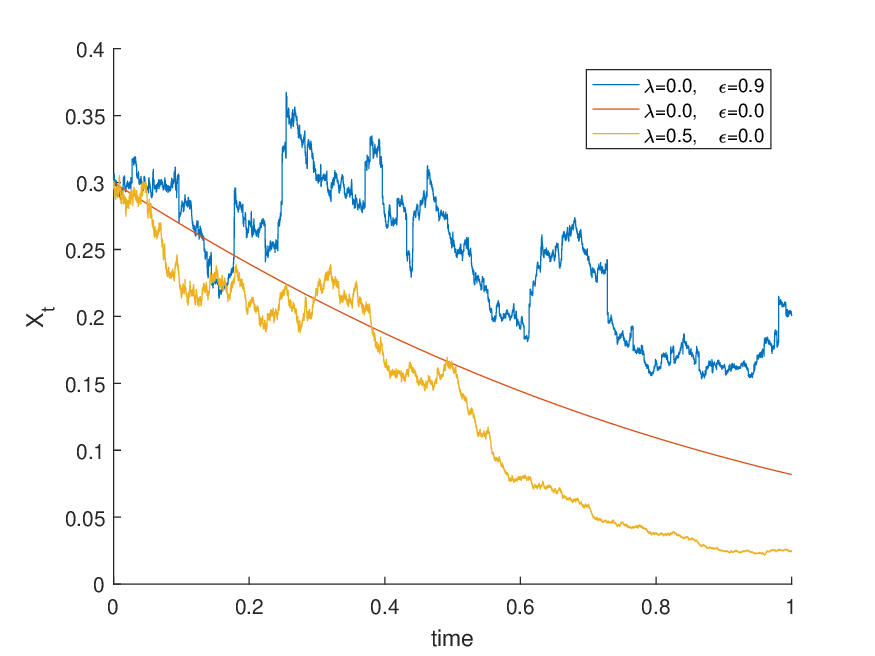}}
	\end{minipage}
\caption{The numerical simulation of the system (\ref{401})  when it is persistent or extinct at different value initial condition $x_0$. (a) Persistent sample paths of model (\ref{401}): the initial condition is {5}.  (b) Extinct sample paths of model (\ref{401}): the initial condition is {0.3}. Parameters $s=1,$\quad $\gamma_2=0.1$,\quad $\gamma_3=2.67$,\quad $\gamma_4=1$,\quad $\alpha=1.5$,\quad $\beta=  0.27< 1$.} \label{path}
\end{figure}

The stochastic single-species system (\ref{401}) with Allee effect, as proved by Jin \cite{7},  has a unique global and positive solution with the initial condition $X_0=x_0$. The jump-diffusion process $X_t$ is adapted and c$\grave{\textup{a}}$dl$\grave{\textup{a}}$g; see Fig. \ref{path}.
Denote the space of c\`{a}dl\`{a}g paths starting at $x_0$ of a solution process $X=\{X_t, t\geq0\}$ of (\ref{401}) by
\begin{equation*}
  \mathcal{D}_{x_0}=\{X :\ \mbox{for any} ~t\geq0,\ \lim_{s\uparrow t}X_s=X_{t-},\  \lim_{s\downarrow t}X_s=X_t\ \mbox{exist}\ \mbox{and} \ X_0=x_0  \}.
\end{equation*}
This space equipped with Skorokhod's $\mathcal{J}_{1}$-topology generated by the metric $\rm{d}_{\mathbb{R}^{+}}$ is a Polish space \cite{13}. For functions $x_{1}, x_{2}\in\mathcal{D}_{x_0}$, define
\begin{eqnarray*}
  {\rm d}_{\mathbb{R}^{+}}(x_{1},x_{2})=\inf\Big\{\varepsilon>0:
     |x_{1}(t)-x_{2}(\bar{\lambda}\, t)|\leq\varepsilon,\ \left|\ln\frac{\arctan(\bar{\lambda}\, t)-\arctan(\bar{\lambda}\, s)}{\arctan(t)-\arctan(s)}\right|\leq\varepsilon,\Big.
     \\
    \left.\mbox{for every}~t, s\geq0\ \mbox{and some}\ \bar{\lambda} \in \Lambda^{\mathbb{R}^{+}}\right\},
\end{eqnarray*}
where $
\Lambda^{\mathbb{R}^{+}}=\{\bar{\lambda} :\mathbb{R}^{+}\rightarrow\mathbb{R};\bar{\lambda}\ \mbox{is injective increasing},\ \lim\limits_{t\rightarrow 0}\bar{\lambda}(t)=0,\ \lim\limits_{t\rightarrow\infty}\bar{\lambda}(t)=\infty \}$.

We consider the corresponding jump-diffusion process $X_{t}(\omega):=\omega(t), t\in[0,T]$ defined on the canonical probability space $({\mathbb{R}}^{[0,T]},\mathcal{B}(\mathbb{R})^{[0,T]},\mathbb{P}_{T})$. Since the paths of $X$ are c\`{a}dl\`{a}g, we identify $X_{t}$  on the space
$(\mathcal{D}_{x_0}^T,\mathcal{B}_{x_0}^T,\mathbb{P})$ instead of $({\mathbb{R}}^{[0,T]},\mathcal{B}(\mathbb{R})^{[0,T]},\mathbb{P}_{T})$, where $\mathcal{D}_{x_0}^T$ is defined similarly as the space $\mathcal{D}_{x_0}$  on the time interval $[0,T]$. The associated Borel $\sigma$-algebra is $\mathcal{B}_{x_0}^T=\mathcal{B}(\mathbb{R})^{[0,T]}\cap\mathcal{D}_{x_0}^T$, and then $(\mathcal{D}_{x_0}^T,\mathcal{B}_{x_0}^T)$ is a separable metric space \cite[Section A.2]{14}. The probability measure $\mathbb{P}$ is generated by $\mathbb{P}(A\cap\mathcal{D}_{x_0}^T):=\mathbb{P}_{T}(A)$ for each $A\in\mathcal{B}(\mathbb{R})^{[0,T]}$. Because every c\`{a}dl\`{a}g function on $[0,T]$ is bounded, we equip $\mathcal{D}_{x_0}^T$
with the uniform norm
\begin{equation*}
\|x\|=\sup_{t\in[0,T]}|x(t)|,\quad x(t)\in\mathcal{D}_{x_0}^T.
\end{equation*}
Hence, $\mathcal{D}_{x_0}^T$ is a Banach space. In order to find the most probable tube of $X_{t}$, we should determine the probability that paths lie within the closed tube
\begin{equation}\label{tube}
K(z,\varepsilon)=\big\{x\in\mathcal{D}_{x_0}^T :\, \parallel x-z\parallel\leq\varepsilon,\,z\in\mathcal{D}_{x_0}^T,\,\varepsilon>0\big\}.
\end{equation}
It is a subset of the space $\mathcal{D}_{x_0}^T$ of c\`{a}dl\`{a}g functions on the interval from $0$ to $T$ containing
a function $z$ together with its $\varepsilon$-neighborhood. Define the measure $\mu_{X}$ on $\mathcal{B}(\mathbb{R})^{[0,T]}$ induced by the solution process $X_t$ for the stochastic nonlinear model (\ref{401}) via
\begin{equation*}
\mu_{X}(B)=\mathbb{P}(\{w: X_{t}(\omega)\in B\}),\quad\text{for}\,\, B\in\mathcal{B}(\mathbb{R})^{[0,T]}.
\end{equation*}
For sufficiently small $\varepsilon>0$, the main contribution of the above probability is given by the measure of the trajectories in the $\varepsilon$-tube of $z\in\mathcal{D}_{x_0}^T$:
\begin{equation}\label{muk}
\mu_{X}(K(z,\varepsilon))=\mathbb{P}(\{w: X_{t}(\omega)\in K(z,\varepsilon)\}),
\end{equation}
where $K(z,\varepsilon)\in\mathcal{B}(\mathbb{R})^{[0,T]}$. As the $\varepsilon$-tube $K(z,\varepsilon)$ depends on the reference path $z$, it is necessary for us to look
for the ``most probable" trajectory $z$ which maximizes the measure $\mu_{X}(K(z,\varepsilon))$ in Eq. \eqref{muk}. When we focus on the differentiable functions $z\in\mathcal{D}_{x_0}^T$, we have the following meaningful definition.

\begin{definition}\label{omf}
Let $0<\varepsilon\ll1$ be given. For a $\varepsilon$-tube surrounding a reference path $z(t)$, the probability of the solution process $X_t,t\in[0,T]$ lying in this tube is estimated by
\begin{equation*}
\mathbb{P}(\parallel X-z\parallel\leq\varepsilon)\varpropto C(\varepsilon)\exp\left\{-\frac{1}{2}\int_{0}^{T}\textup{OM}(\dot{z},z)dt\right\},
\end{equation*}
 where the integrand $\textup{OM}(\dot{z},z)$ is called Onsager-Machulup function and $\varpropto$ denotes the equivalence relation for $\varepsilon$ small enough. The intergral $\int_{0}^{T}\textup{OM}(\dot{z},z)dt$ is the Onsager-Machulup functional.
\end{definition}

\begin{remark}
The Onsager-Machulup function is similar to the Lagrangian function of a dynamical system in classical mechanics, and the OM functional would correspond to the action functional. In particular, for an SDE with pure jump L\'evy
noise, Definition \ref{omf} is still applicable, and the minimizer of the OM functional $\int_{0}^{T}\textup{OM}(\dot{z},z)dt$ gives the most probable path for this non-Gaussian stochastic system. Moreover, the minimizer $z$ may be chosen from a more general function space.
\end{remark}

Our main result about the expression of the OM function for a jump-diffusion process is clearly presented in the basic theorem.
\begin{theorem}\label{theoremOM}
For the stochastic nonlinear system \eqref{401} with the jump measure satisfying $\int_\mathbb{Y}\epsilon(y)\nu_{\alpha}(dy)<\infty$, the Onsager-Machlup function \cite{15} is characterized, up to an additive constant, by:
{\small\begin{align*}
 \textup{OM}(\dot{z},z)=&\left[\frac{\dot{z}-z\left(s-\gamma_2\,z-\frac{\gamma_3}{\gamma_3\,\gamma_4\,z+1}\right)}{\lambda z}\right]^{2}+s-2\gamma_2z-\frac{\gamma_3}{(\gamma_3\,\gamma_4\,z+1)^{2}}\\
 &+2\frac{\dot{z}-z\left(s-\gamma_2\,z-\frac{\gamma_3}{\gamma_3\,\gamma_4\,z+1}\right)}{\lambda^2 z}\int_\mathbb{Y}\epsilon(y)\nu_{\alpha}(dy),
\end{align*}}
where $z\in\mathcal{D}_{x_0}^T$ is a differentiable function. The contribution of pure jump L\'evy noise to the OM function is the third term. When the jump measure is absent, we cover the OM function for the case of diffusion. In terms of OM function, the measure of tube $K(z,\varepsilon)$ defined in \eqref{tube} can be approximated as follows:
\begin{equation*}
\mu_{X}(K(z,\varepsilon))\varpropto \mu_{Y^{c}}(K(0,\varepsilon))\exp\left\{-\frac{1}{2}\int_{0}^{T}\textup{OM}(\dot{z},z)dt\right\},
\end{equation*}
where  $Y_t^{c}$ is defined
by
\begin{equation*}
dY_t^{c}=Y_{t}^{c}\left(\lambda\,Y_t^{c}dB_t+\int_\mathbb{Y}\epsilon(y)\tilde{N}(dt,dy)\right), \quad t\in[0,T].
\end{equation*}
\end{theorem}
The proof of Theorem \ref{theoremOM} is given in \cite[Theorem 4.1]{5}.

In Gaussian noise case ($\epsilon(y)=0)$,  the stochastic single-species model (\ref{401}) becomes

\begin{equation}\label{sdeG}
dX_t=X_{t}\left[\left(s-\gamma_2\,X_{t}-\frac{\gamma_3}{\gamma_3\,\gamma_4\,X_{t}+1}\right)dt + \lambda dB_t\right], \quad t\geq 0,\quad X_0=x_0.
\end{equation}
We apply Lamperti transforms for solving SDE driven by multiplicative noise \cite[Example 6.48]{1}. This method allows us to transform the multiplicative noise into additive noise. Because numerically solving an additive-noise SDE is usually easier than solving a multiplicative-noise SDE as in Eq. (\ref{sdeG}).

Assume $g\in C^2(\mathbb{R})$ and define  $Y_t=g=\ln(X_t)$.  Then the new SDE has the following form:
\begin{equation}\label{Ito-F2}
 dY_t=G(Y_t)dt+\lambda\,dB_t,
\end{equation}
where
$$G(Y_t)=\left (h(X_t)-\frac{\lambda^2}{2} \right) \Big|_{\,X_t\, = \,\exp (Y_t)},$$ and  $$h(X_t)=s-\gamma_2\,X_t-\frac{\gamma_3}{\gamma_3\,\gamma_4\,X_t+1}.$$

Since the most probable transition path for a stochastic single-species model is the minimizer of the Onsager-Machlup action functional, denoted by $Z_m$, it can be obtained from the following least action principle
\begin{align*}
  \delta  \int_{0}^{T}\textup{OM}(\dot{z},z)dt=0,
\end{align*}
where the integrand function (Onsager-Machlup function) \cite{6} is given by
\begin{equation}\label{OMsde}
\textup{OM}(\dot{z},z)=\left( \frac{G(z)-\dot{z}}{\lambda}\right)^2 +\dot{G}(z).
\end{equation}
Thus Eq. (\ref{OMsde}) satisfies the following Euler-Lagrange equation
\begin{align}\label{EL}
   \frac{d}{dt}\frac{\partial \textup{OM}(\dot{z},z)}{\partial \dot{z}} =\frac{\partial \textup{OM}(\dot{z},z)}{\partial z}.
\end{align}
The most probable transition pathway $Z_m(t)$ of system (\ref{Ito-F2}) is characterized by
\begin{align}\label{Zm}
 &\ddot{Z}_{m}(t)=  \frac{\lambda^2}{2}\ddot{G}(Z_m)+\dot{G}(Z_m)\,G(Z_m),\qquad 0 < t < T,\nonumber\\
 &Z_m(0)=X_1,\,\,\, Z_m(T)=X_3.
\end{align}
To solve two-point boundary value problem in Eq. (\ref{Zm}), we apply the shooting method depicted in Ref. \cite{sh}.

\subsection{Most probable phase portraits}\label{Mppp}
As for the solution of the Fokker-Planck equation, the probability density function $p(X,t)$ is a surface in the $(X,t,p)$-space. For a given time $t$, the maximizer $X_m(t)$ for $p(X,t)$  (i.e., $X_m(t)=\mathrm{max}_{X \in (0,\infty)}p(X,t)$) shows the most probable (i.e., maximal likely) location of this orbit at time $t$. The orbit traced out by $X_m(t)$ is called a most probable orbit starting at $x_0$. Thus, the deterministic orbit $X_m(t)$ follows the top ridge of the surface in the $(X,t,p)$-space as time goes on.
\subsubsection{Non-local Fokker-Plank equation}
{The Fokker-Planck equation describes the time evolution of the probability density function, but it can be solved analytically only in special cases. We are interested in the steady-state probability distribution (equilibrium distribution), and want to express the stationary solution of the non-local Fokker-Planck equation. This makes the estimate of the most probable phase portrait possible in L\'evy noise case numerically and algorithmically.}

Let $f:\mathbb{R}\rightarrow\mathbb{R}$ be a smooth function. Suppose that the solution $X_{t}$ of system (\ref{401}) has a conditional probability density $p(X,t|x_{0},0)$. For convenience, we drop the initial condition and simply denote it by $p(X,t)$. On one hand,
\begin{equation*}
\mathbb{E}f(X_{t})=\int_{\mathbb{R}}f(X)p(X,t)dX,
\end{equation*}
and thus
\begin{equation*}
\frac{d}{dt}\mathbb{E}f(X_{t})=\int_{\mathbb{R}}f(X)\frac{\partial}{\partial t}p(X,t)dX.
\end{equation*}
On the other hand, by virtue of It\^o's formula,
\begin{align}\nonumber
df(X_{t})=&X_{t-}\left(s-\gamma_2\,X_{t-}-\frac{\gamma_3}{\gamma_3\,\gamma_4\,X_{t-}+1}\right)f'(X_{t})dt\\ \label{df}
          &+\int_{\mathbb{Y}}\big(f(X_{t}+\epsilon(y)X_{t-})-f(X_{t})-\epsilon(y)X_{t-}f'(X_{t})\big)\nu_{\alpha}(dy)dt.
\end{align}
Taking expectation on both sides of (\ref{df}), we gain
\begin{align}\nonumber
d\mathbb{E}f(X_{t})=&\mathbb{E}\Big[X_{t-}\left(s-\gamma_2\,X_{t-}-\frac{\gamma_3}{\gamma_3\,\gamma_4\,X_{t-}+1}\right)f'(X_{t})dt\\\label{ef1}
  &+\int_{\mathbb{Y}}\big(f(X_{t}+\epsilon(y)X_{t-})-f(X_{t})-\epsilon(y)X_{t-}f'(X_{t})\big)\nu_{\alpha}(dy)dt\Big].
\end{align}
Noting that the infinitesimal generator of the solution $X_{t}$ for system (\ref{401}) is
\begin{align*}
 Ap(X,t)=&X\left(s-\gamma_2\,X-\frac{\gamma_3}{\gamma_3\,\gamma_4\,X+1}\right)\partial_Xp(X,t)\\
         &+\int_{\mathbb{Y}}\big(f(X+\epsilon(y)X)-f(X)-\epsilon(y)X\partial_Xp(X,t)\big)\nu_{\alpha}(dy).
\end{align*}
The equation (\ref{ef1}) is rewritten as
\begin{align}\nonumber
\frac{d}{dt}\mathbb{E}f(X_{t})=
&\mathbb{E}\Big[X_{t-}\left(s-\gamma_2\,X_{t-}-\frac{\gamma_3}{\gamma_3\,\gamma_4\,X_{t-}+1}\right)f'(X_{t})\\ \nonumber
&+\int_{\mathbb{Y}}\big(f(X_{t}+\epsilon(y)X_{t-})-f(X_{t})-\epsilon(y)X_{t-}f'(X_{t})\big)\nu_{\alpha}(dy)\Big]\\ \nonumber
=&\int_{\mathbb{R}}\Big[X\left(s-\gamma_2\,X-\frac{\gamma_3}{\gamma_3\,\gamma_4\,X+1}\right)f'(X)\\ \label{fpen-1}
&+\int_{\mathbb{Y}}\big(f(X+\epsilon(y)X)-f(X)-\epsilon(y)X\,f'(X)\big)\nu_{\alpha}(dy)\Big]p(X,t)dX.
\end{align}
As a result, the Fokker-Planck equation for the stochastic nonlinear system (\ref{401}) of the solution process $X = \{X_t, t\geq 0\}$ with initial condition $p(X,0)=\sqrt {\frac{40}{\pi}}e^{-40\left(X-x_0\right)^2}$ is
{\small\begin{align}\label{FK}
\partial_tp(X,t)=&-\left(s-2\gamma_2X-\frac{\gamma_3}{(\gamma_3\,\gamma_4\,X+1)^{2}}\right)p(X,t)-X\left(s-\gamma_2\,X -\frac{\gamma_3}{\gamma_3\,\gamma_4\,X+1}\right)\partial_Xp(X,t)\nonumber\\
  &+\int_{\mathbb{Y}}\big(f(X+\epsilon(y)X)-f(X)-\epsilon(y)X\,f'(X)\big)\nu_{\alpha}(dy)p(X,t).
\end{align}}

To simulate the
non-local Fokker-Planck equation (\ref{FK}), we apply a numerical finite difference method given in Gao et al. \cite{Gao}.

If the L\'evy motion is replaced by Brownian motion, then the local Fokker-Planck equation has the following form:
\begin{equation}\label{BMFK}
\partial_tp(X,t)=-\partial_{X}\left[X\left(s-\gamma_2\,X -\frac{\gamma_3}{\gamma_3\,\gamma_4\,X+1}\right)p(X,t)\right]+\frac{\lambda^2}{2}\partial_{XX}[X^2\,p(X,t)].
\end{equation}

\begin{figure}[H]
	\begin{minipage}{0.48\linewidth}
		\leftline{(a)}
		\centerline{\includegraphics[height = 6cm, width = 8.6cm]{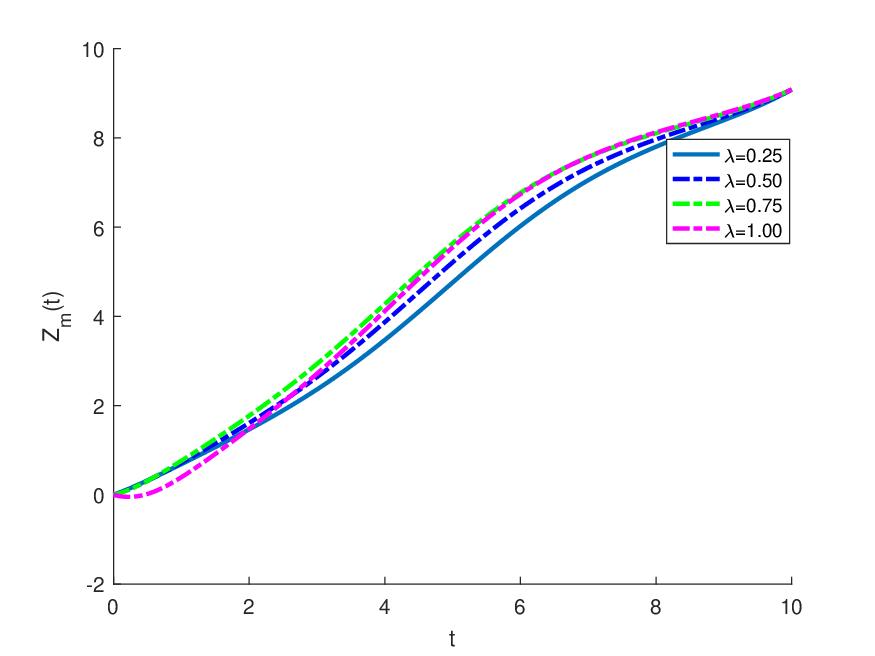}}
	\end{minipage}
	\hfill
	\begin{minipage}{0.48\linewidth}
		\leftline{(b)}
		\centerline{\includegraphics[height = 6cm, width = 8.6cm]{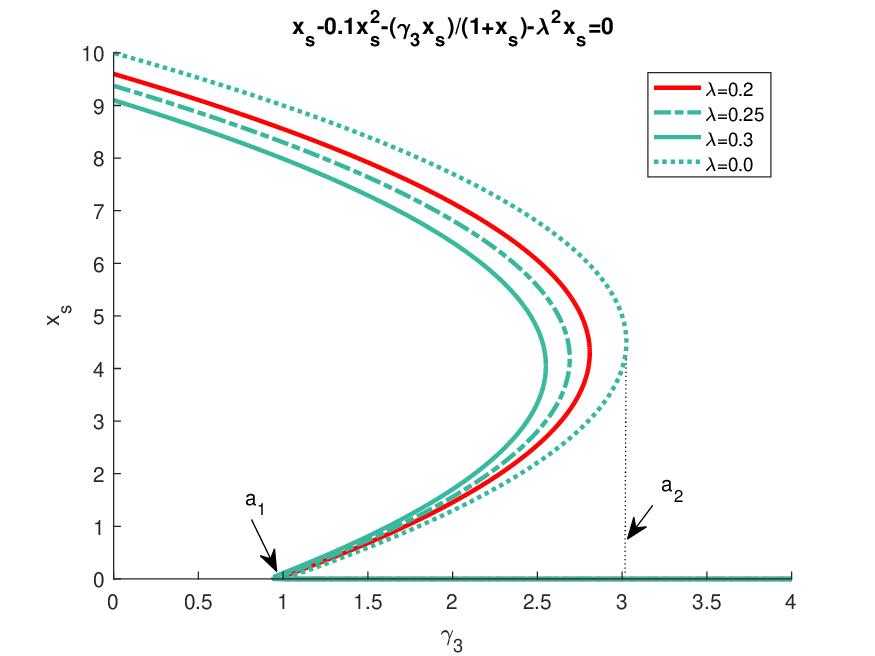}}
	\end{minipage}
\caption{(a) Most probable transition pathways $Z_m(t)$ starting at the extinction state $X_1=0$ and ending at the upper equilibrium stable state $X_3=9.0846$ under white noise with respect to the different values of $\lambda$. (b) The most probable steady state $x_s$ versus the attack rate $\gamma_3$ for different values of Gaussian noise intensity with $\lambda$ term.}\label{OMBM}
\end{figure}

The stationary probability density function $p_s(X)$ of Eq. (\ref{BMFK}) can be solved by
\begin{equation}\label{FKnBM}
0=-\partial_{X}\left[X\left(s-\gamma_2\,X -\frac{\gamma_3}{\gamma_3\,\gamma_4\,X+1}\right)p_s(X)\right]+\frac{\lambda^2}{2}\partial_{XX}[X^2\,p_s(X)],
\end{equation}
or equivalently,
\begin{align}\label{FM}
&0=-\left[X\left(s-\gamma_2\,X -\frac{\gamma_3}{\gamma_3\,\gamma_4\,X+1}\right)p_s(X)\right]+\frac{\lambda^2}{2}\partial_{X}[X^2\,p_s(X)],\nonumber\\
 &\Longrightarrow\quad 0=\left[X\left(s-\gamma_2\,X -\frac{\gamma_3}{\gamma_3\,\gamma_4\,X+1}\right)-\lambda^2\,X\right]p_s(X)-\frac{\lambda^2\,X^2}{2}\partial_{X}\,p_s(X).
  \end{align}
Due to the complexity of stationary solution, we take the extrema of the stationary probability density  function (spdf) located at $x_s$ directly. In other words,  the spdf satisfies $\partial_{X}(\,p_s(x_s))=0$. Since $p_s(x_s)\neq0,$
Eq. (\ref{FM}) becomes
\begin{equation}\label{xs}
X\left(s-\gamma_2\,X -\frac{\gamma_3}{\gamma_3\,\gamma_4\,X+1}\right)-\lambda^2\,X=0.
\end{equation}
Eq. (\ref{xs}) is completely different from the equilibrium state of the deterministic model (\ref{42}) because of the presence of noise with $\lambda$ term. The numerical solution of Eq. (\ref{xs}) is plotted in Fig. \ref{OMBM}(b).

\section{Numerical results and its biological implications}\label{numerical}

In order to make the readers understand our results much better, we perform some numerical simulations to illustrate our theoretical results. Based on the finite difference method \cite{Gao}, the numerical simulations are very useful in the study of real examples of population.
In the present section, we define the bifurcation time as the time between the changes in number of maximal likely equilibrium states. It is a time scale for the birth of a new most probable stable equilibrium state. We also show the intervals where there exist one or two maximal likely stable equilibrium states, the value of the equilibrium states, and the point where the number of metastable states of the stochastic single-species model \eqref{401} varies.
Since the numerical solutions of a model depend on the values of all its deterministic parameters and noise intensities. Here, we discuss the effect of the parameters in Table \ref{P} on the investigated system \eqref{401}. For simplicity, we choose four maximum likely pathways together with the initial conditions selected in different intervals.

While we plot the above figures, we  fix  the deterministic parameters $s=1$,\, $\gamma_2=0.1$,\, $\gamma_3=2.67,$\, $\gamma_4=1/\gamma_3,$  the noise intensity  $\epsilon=0.5,$ and the stability index $\alpha=1.5$.

The potential function denoted by $U(X)$  in Fig. \ref{de-pot-bif}(a) has two stable steady states $X_1$ and $X_3$, and an unstable steady state $X_2$ for $\beta <1.$ This function  has a maximum value at the unstable equilibrium solution $X_2$.  At the stable fixed points $X_1$ and $X_3$, the potential function attains its minima. For the value of  $\beta > 1,$ the nonlinear system (\ref{42}) has only one equilibrium point, which is the trivial point $X_1=0$.

In Fig. \ref{de-pot-bif}(b), we sketch the equilibrium states versus attack rate $\gamma_3$. For $\beta < 1$, there exist two stable equilibrium states $X_1$ and $X_3$ and one unstable equilibrium state $X_2.$ While $\beta >1$, $X_1 = 0$ is the unique equilibrium state that is stable. Thus, the parameter $\beta=1$ is the bifurcation parameter value.

 The distance between the unstable equilibrium $X_2$ and the stable fixed point $X_3$ becomes very small when $\beta$ approaches to 1. This indicates that the expected time to extinction may be too short, as clarified in Fig. \ref{de-pot-bif}(b).

{Fig. \ref{path} displays the numerical simulation of the stochastic single-species model (\ref{401}) with Allee effect  when it is persistent or extinct at different value of initial condition $x_0$. This figure proves  that the solutions of the stochastic nonlinear system \eqref{401} are positive, and extinction species occurs when the initial condition is less than the value of $X_2$, as demonstrated in Fig. \ref{path}(b). While the initial condition is greater than the value of $X_2$, there is stochastic persistence.}

In Fig. \ref{OMBM}(a), we depict the most probable transition pathways $Z_m(t)$ of system (\ref{Ito-F2}) for different values of $\lambda$. {This figure tells us that as time grows, the most probable paths $Z_m(t)$ increase to the stable state $X_3$ quickly, and remain at a nearly constant level, then approach to the high stable equilibrium state}. Fig. \ref{OMBM}(b) demonstrates the  curves for the most probable steady state $x_s$ of the stochastic single-species model \eqref{401} with $\epsilon(y)=0$ driven by Gaussian noise at different values of the noise intensity $\lambda.$  The steady state curves exhibit a bi-stability in the interval $(a_1, a_2)$. For $\gamma_3 > a_1$, the stable steady state stays  at the extinction state. {While for $\gamma_3 <a_1$, it is located at the stable equilibrium state. Because of the presence of Gaussian noise with $\lambda$ term,} the numerical  result in Fig. \ref{OMBM}(b) is completely different from the numerical result in Fig. \ref{de-pot-bif}(b).

\begin{figure}[H]
\centering
{\includegraphics[width=0.7\textwidth]{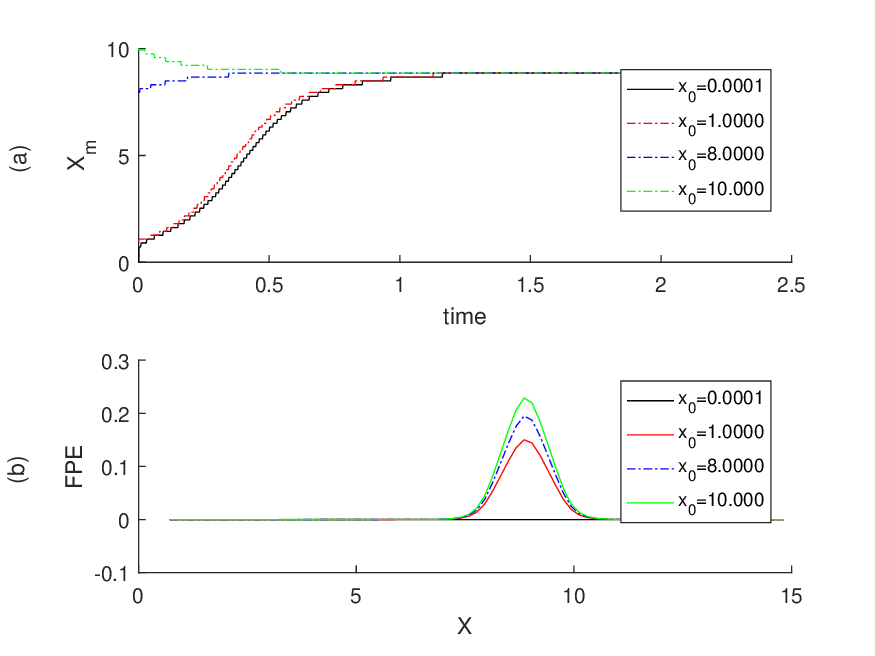}}
\caption{(a) Most probable orbits and most probable equilibrium states for system (\ref{401}) with equilibrium state $X_m$ between 9 and 10. (b) The solution of Fokker-Planck equation of model (\ref{401}). The stationary density function of the FPE has maximum value at the equilibrium state $X_m$. The other parameters are fixed: $s=1,$\quad $\gamma_2=0.1$,\quad $\gamma_3=2.67$,\quad $\gamma_4=1$, $\alpha=1.5$, $\epsilon=0.5,$ $\beta= 0.27< 1$,\,\,$\lambda=0$ and $x_0\in (0, 10]$. The bifurcation time is 1.13.}\label{fpe-mppN1}
\end{figure}

Fig. \ref{fpe-mppN1}(b) draws the MPPP for different initial values $x_0$. From this figure, we observe that the maximum value of the stationary density function $p(X,t)$ is located at the maximum likely stable state $X_m(t)=9.0846$ with the initial condition $p(X,0)=\sqrt {\frac{40}{\pi}}e^{-40\left(X-x_0\right)^2}$. As the initial condition $x_0$ increases, it raises the peak point of the stationary density function  $p(X,t)$. This shows that the extinction of the species may not happen, and the high peak occurs at the maximum likely stable state $X_m(t)$.
\begin{figure}[H]
\begin{center}
  \begin{minipage}{2.1in}
\leftline{(a)}
\includegraphics[width=2.1in]{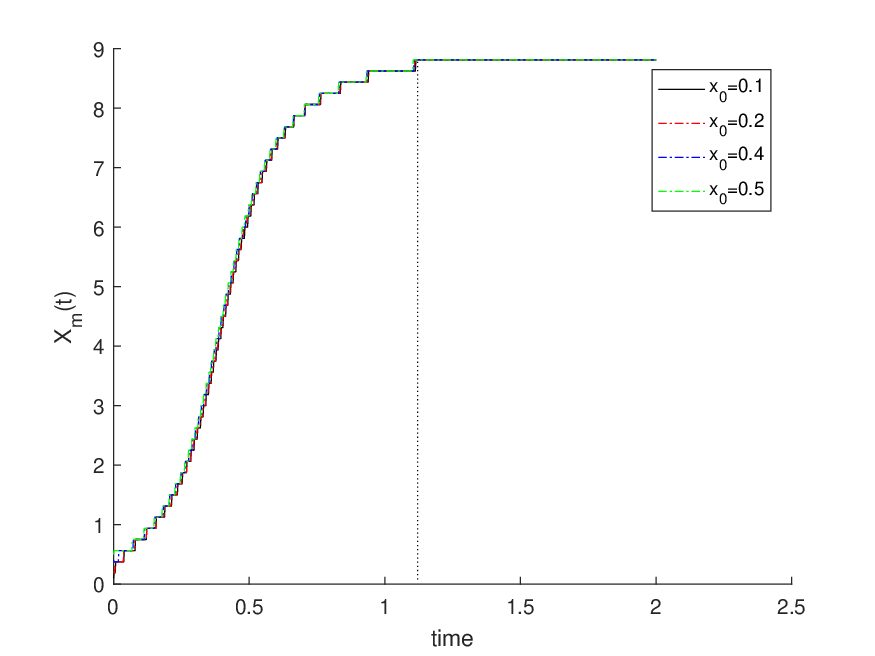}
\end{minipage}
\hfill
\begin{minipage}{2.1in}
\leftline{(b)}
\includegraphics[width=2.1in]{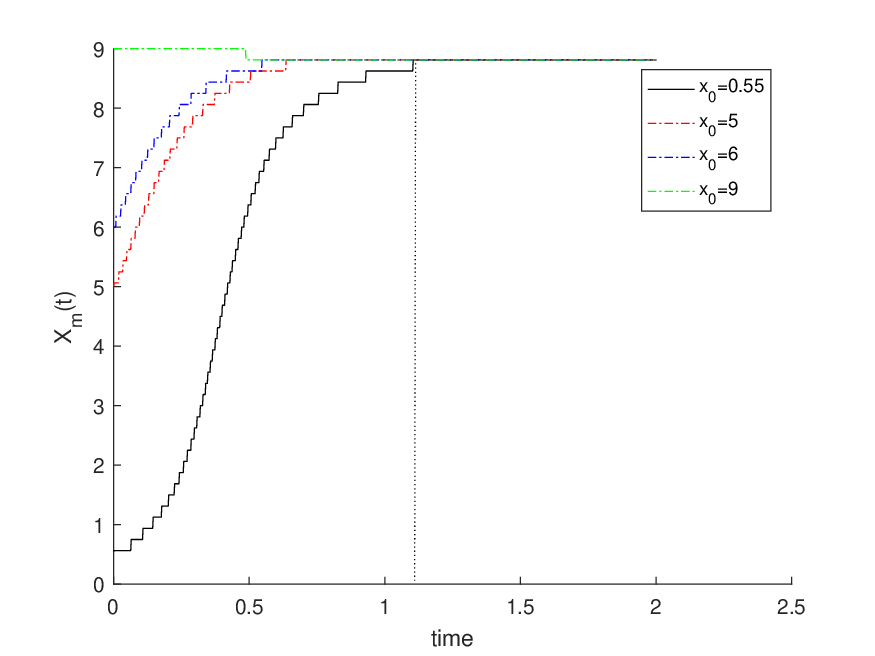}
\end{minipage}
\hfill
  \begin{minipage}{2.1in}
\leftline{(c)}
\includegraphics[width=2.1in]{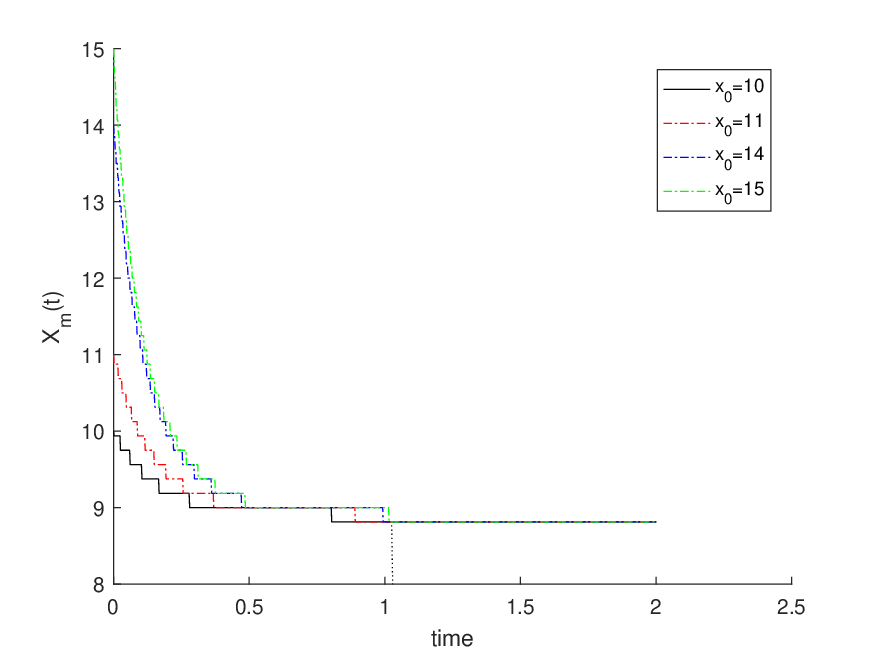}
\end{minipage}
\caption{Most probable orbits and most probable equilibrium states for stochastic nonlinear system (\ref{401}). (a) When the initial condition $x_0$ is less than the unstable equilibrium state $X_2$, i.e., $x_0 < X_2$. (b) When the initial condition is between $X_2$ and $X_3$, i.e., $X_2 < x_0 < X_3$. (c) When the initial condition is greater than $X_3$, i.e., $x_0>X_3$. Parameters $s=1,$\quad $\gamma_2=0.1$,\quad $\gamma_3=2.67$,\quad $\gamma_4=1$,\quad $\alpha=1.5$,\quad $\epsilon=0.5,$ \quad $\beta=0.27< 1$, $\lambda=0$, and bifurcation time at 1.13 (dot vertical line).}\label{max}
\end{center}
\end{figure}

The most probable trajectories of the stochastic single-species model (\ref{401}) with Allee effect are plotted graphically in figures \ref{max} and \ref{fpe-mppN1} (a). Here the values of the noise intensities are set up as  $\lambda=0$ and $\epsilon=0.5,$ respectively. We choose the stability index $\alpha=1.5,$ and the interval $D=(0,15)$. These figures evolve as the initial value $x_0$ changes, and they tell us that  the  maximal likely equilibrium state (maximizer) $X_m(t)$ lies between 9 and 10  at the bifurcation time  1.13. In other words, the maximizer in high concentration is between 9 and 10, it  is  different from the deterministic equilibrium stable solution $X_3=9.0846$ due to the effect of external noises. As seen in Fig. \ref{fpe-mppN1}(a), the most probable growth state is attracted to the maximal likely equilibrium state in the extinction state, and then it leads to the maximal likely equilibrium state in the high concentration as time moves forward.

\section{Conclusion }\label{conclusion}
In the present work, we have studied the Onsager-Machlup functional and  most probable phase portraits for the stochastic growth model \eqref{401} for single-species population with strong Allee effects driven by L\'evy noise. We have focused on the effect of different values of the initial condition on the MPPP of the nonlinear dynamical system. Small disturbances may cause a transition between the extinction stable state $X_1$ and the upper equilibrium state $X_3$, thus we have used a deterministic quantity, namely the maximal likely trajectory to analyze the transition phenomena in a jump stochastic environment.

In order to find the most likely pathways in transition phenomena, we have calculated the most probable paths of stochastic differential equation in (\ref{401}) using the stationary density function of the non-local Fokker-Planck equation associated with a non-local partial differential equation. We have investigated the impact of the deterministic  parameters, noise intensities and domain size  on the FPE. We also have studied the dependence of the probability density on the initial condition $x_0$. Our finding has displayed that the maximum of the stationary density function is located at the most probable stable equilibrium state $X_m$.

In conclusion, the most probable path has been used as an indicator that helps the researcher to understand the stochastic dynamics of the single-species model \eqref{401} based on the evolution of the probability density function over time.

\bigskip
\noindent\textbf{Data Availability}

Numerical algorithms source code that support the findings of this study are openly
available in GitHub, Ref. \cite{Y21}.

\bigskip
\noindent\textbf{Acknowledgements}

The authors are happy to thank Professor Jinqiao Duan for fruitful discussions on stochastic dynamical systems. The authors acknowledge support from the NSFC grant 12001213.

\bigskip

\noindent\textbf{Conflict of Interest}

The authors declare that they have no conflict of interest.

\end{document}